\documentclass[11pt]{article}
\usepackage{amsmath,amsfonts,latexsym,amssymb,amscd, graphicx,pictexwd, mathrsfs, dsfont, color,mathtools}
\newcommand{\F}{{\mathbf F}}
\newcommand{\B}{{\mathbf B}}

\newcommand{\LL}{{\mathbf L}}

\renewcommand{\P}{{\mathbb P}}
\newcommand{\E}{{\mathbb E}\,}
\newcommand{\bL}{{\mathbb L}\,}

\newcommand{\R}{{\mathbb R}}

\newcommand{\Q}{{\mathbb Q}}

\newcommand{\KK}{{\mathit K}}
\newcommand{\II}{{\mathit I}}

\newcommand{\ZZ}{{\mathbb Z}}

\newcommand{\cA}{{\mathcal A}}
\newcommand{\cL}{{\mathcal L}}

\newcommand{\cB}{{\mathcal B}}

\newcommand{\cP}{{\mathcal P}}
\newcommand{\cC}{{\mathcal C}}

\newcommand{\fh}{{\mathfrak h}}
\newcommand{\fb}{{\mathfrak b}}

\newcommand{\ff}{{\mathfrak f}}
\newcommand{\fd}{{\mathfrak d}}

\newcommand{\sfUC}{{\mathsf{UC}}}

\newcommand{\argmax}{\mathop{\rm arg\,max}}

\newtheorem{thm}{Theorem}

\newtheorem{lem}{Lemma}[section]
\newtheorem{prop}{Proposition}
\newtheorem{Definition}{Definition}
\newtheorem{rem}{Remark}

\newtheorem{Assump}{Assumption}

\numberwithin{equation}{section}

\begin{document}
\title{Brownian Aspects of the KPZ Fixed Point}
\author{Leandro P. R. Pimentel}
\date{\today}
\maketitle

\begin{abstract}  
The Kardar-Parisi-Zhang (KPZ) fixed point is a Markov process $\left(\fh_t\,,\,t\geq 0\right)$ that is conjectured to be at the core of the KPZ universality class \cite{CoQuRe}. In this article we study two aspects the KPZ fixed point that share the same Brownian limiting behaviour: the local space regularity and the long time evolution. Most of the results that we will present here were obtained by either applying explicit formulas for the transition probabilities \cite{MaQuRe} or applying the coupling method to discrete approximations \cite{Pi,Pi2}. Instead we  will use the variational description of the KPZ fixed point \cite{DaOrVi}, allowing us the possibility of running the process starting from different initial data (basic coupling), to prove directly the aforementioned limiting behaviours.  
\end{abstract}

\section{Introduction}
 
The universality class concept is an artifact of modern statistical mechanics that systemizes the  idea that there are but a few important characteristics that determine the scaling behaviour of a stochastic model.  In $d+1$ stochastic growth models the object of interest is a height function $h(x,t)$ over the $d$-dimensional substrate $x\in\R^d$ at time $t\geq 0$, whose evolution is described by a random mechanism. For fairly general models one has a deterministic macroscopic shape for the height function and its fluctuations, under proper space and time scaling, are expected to be characterized by a universal distribution. A well known example is given by the random deposition growth model, where  blocks are pilled in columns (indexed by $x\in\ZZ$) according to independent Poisson processes. The existence of a macroscopic shape follows from the law of large numbers and, due to the classical central limit theorem, the height function at $x\in\R$ has Gaussian fluctuations that are independent in space. In 1986 \cite{KPZ}, Kardar, Parisi and Zhang (KPZ) proposed a stochastic partial differential equation (the KPZ equation) for a growth model where a non-linear local slope dependent rate is added to a diffusion equation with additive noise: $\partial_t h=\frac{1}{2}(\partial_x h)^2+\partial^2_x h+\xi $. In opposition to the previous random deposition growth model, they predicted that for $d=1$ the solution of the KPZ equation at time $nt$ has fluctuations of order $n^{1/3}$, and on a scale of $n^{2/3}$ that non-trivial spatial correlation is achieved (KPZ scaling exponents). Since then it is expected that $1+1$ interface growth models that exhibit a similar KPZ growth mechanism would satisfy 
$$h(an^{2/3}x,nt)\sim bnt+cn^{1/3}\fh_t(x)\,,$$
for some constants $a,b,c\in \R\setminus\{0\}$ that might depend on the microscopic distributional details of the model, but where $\fh_t(x)$ is a universal space-time process called the KPZ fixed point \cite{CoQuRe}. Illustrations of natural phenomena within the KPZ universality class include turbulent liquid crystals, bacteria colony growth and paper wetting \cite{Ta}. For a more complete introduction to the KPZ equation and universality class, and its relation with other discrete growth models in statistical physics, the author address to \cite{Co}.   
\newline

After \cite{KPZ}, the study of KPZ fluctuations became a notorious subject in the literature of physics and mathematics and, in the late nineties, a breakthrough was presented by Baik, Deift and Johansson \cite{BaDeJo, Jo1}. By applying an exact formula (in terms of a Toeplitz determinant) for the Hammersley last-passage percolation growth model with narrow wedge initial profile, and then by analysing asymptotics of the resulting expressions, they were able to prove convergence of shape fluctuations at $x=0$ to the Tracy-Widom (GUE) distribution. In the past twenty years there has been a significant amount of improvements of the theory and the exact statistics for some special initial conditins, resulting in different types of limiting distributions, were computed using connections with integrable probability \cite{AmCoQu,BoFePrSa,Jo2,PrSp}. Recently, a unifying approach was developed by Matetski, Quastel and Remenik \cite{MaQuRe} in the TASEP\footnote{Totally asymmetric simple exclusion process.} growth model context that conducted to the formal construction of the Markov process $(\fh_t\,,\,t\geq 0)$ and the explicit computation of the transition probabilities.
\newline

Alongside the rich structure of integrable probability, the study of the KPZ universality class was also developed by techniques based on the graphical representation of an interacting particle system due to Harris \cite{Ha}. There are many conveniences of this approach, also known as the coupling method, comprising the possibility of running the process starting from different initial data on the same probability space. In the seminal paper by Cator and Groeneboom \cite{CaGr}, the authors applied the coupling method to derive the KPZ scaling exponents ($1/3$ and $2/3$) for the Hammersley last-passage percolation growth model. This method was further developed in the TASEP context by Bal\'azs, Cator and Sepp\"al\"ainen \cite{BaCaSe}, and became a successful tool to analise fluctuations of  models \cite{BaQuSe,BaSe,Se} lying within the KPZ universality class, and local properties of different types of of Airy processes \cite{CaPi,FeOc,Pi}. Related to that, there has been considerable developments in describing the space-time structure of the KPZ fixed point in terms of a variational formula \cite{CoLiWa,FeOc} that relies on the existence and uniqueness of a two-dimensional random scalar field, called the Airy sheet \cite{CoQuRe,DaOrVi,MaQuRe}. In analogy with Harris graphical representation, this variational formula allows one to run simultaneously the process starting from different initial data on the same probability space (basic coupling). Thereby, it seems natural to expect that particle systems techniques that were applied to discrete approximations of the KPZ fixed point \cite{Pi,Pi2} can be developed in the continuous space-time context itself. In the course of this article we prove Brownian behaviour of the KPZ fixed point (Theorem \ref{Holder_1}, Theorem \ref{AirySheet} and Theorem \ref{Coupling} in Section \ref{results}) by using soft arguments based on geometrical aspects of the variational representation formalized by Dauvergne, Ortmann and Virag \cite{DaOrVi}.   
\newline

\noindent\paragraph{\bf Acknowledgements} This research was supported in part by the International Centre for Theoretical Sciences (ICTS) during a visit for participating in the program - Universality in random structures: Interfaces, Matrices, Sandpiles (Code: ICTS/urs2019/01), and also by the CNPQ grants 421383/2016-0 and 302830/2016-2, and by the FAPERJ grant E-26/203.048/2016.   

\section{Main Results}\label{results}

\subsection{The Airy Sheet and the Directed Landscape}
The construction of the directed landscape is based on the existence and uniqueness of the so called Airy Sheet, which in turn is defined through a last-passage percolation model over the parabolic Airy line ensemble \cite{DaOrVi}. For a sequence of differentiable functions $\F=(\dots,\F_{-1},\F_0,\F_1,\dots)$ with domain $\R$, and coordinates $x\leq y$ and $n\leq m$, define the last-passage percolation time
$$\F\left((x,m)\to(y,n)\right):=\sup_{\pi}\int_x^y\F'_{\pi(t)}(t)dt \,,$$   
where the supremum is over nonincreasing functions $\pi:[x,y]\to\ZZ$ with $\pi(x)=m$ and $\pi(y)=n$. Notice that, for such paths, the integral is just the sum of the increments of $\F$ (over each line), so the same can be defined for continuous $\F$. An important exemple is given by setting $\F\equiv\B$ a sequence of independent standard two-sided Brownian motions (Brownian last-passage percolation). In the literature of last-passage percolation it is normally considered maximization over   nondecreasing paths instead, but to accommodate the natural order of the Airy line ensemble from top to bottom (as below), Dauvergne, Ortmann and Virag \cite{DaOrVi} defined it for nonincreasing paths.           
\newline

The parabolic Airy line ensemble \cite{PrSp} is a random sequence of ordered real functions $\LL_1>\LL_2>\dots$ with domain $\R$. The function   $\LL_n(x)+x^2$ is stationary for all $n\geq 1$, and the top line $\LL_1(x)+x^2$ is known as the Airy$_2$ process and represents the limit fluctuations of some integrable last-passage percolation models, including the Brownian one.     

\begin{Definition}
The stationary Airy sheet is a random continuous function $\cA:\R^2\to\R$ such that:
\begin{itemize}
\item $\cA\stackrel{dist.}{=}T_{(z,w)}\cA$ for all $(z,w)\in\R^2$, where $T_{(z,w)}\ff(x,y):=\ff(x+z,y+w)$.
\item $\cA$ can be coupled with the parabolic Airy line ensemble so that 
$$\left(\cA(0,x)\,,\,x\in\R\right) \stackrel{dist.}{=}\left(\LL_1(x)+x^2\,,\,x\in\R\right)\,,$$ 
and for all $(x,y,z)\in\Q^+\times \Q^2$ almost surely there exists a random $K_{x,y,z}$ such that for all $k\geq  K_{x,y,z}$ we have 
$$\LL\left((-\sqrt{k/2x},k)\to (z,1)\right) - \LL\left((-\sqrt{k/2x},k)\to (y,1)\right)=\cL(x,z)-\cL(x,y)\,,$$
where 
\begin{equation}\label{PaAiSh}
\cL(x,y):=\cA(x,y)-(x-y)^2\,
\end{equation}
is called the parabolic Airy sheet.
\end{itemize}
\end{Definition}

In \cite{DaOrVi} the authors have a similar definition for the parabolic Airy sheet $\cL(x,y)$ (see Definition 1.2\cite{DaOrVi} and notice that they used different notation to represent the parabolic Airy line ensemble and the parabolic Airy sheet), but it follows from their results (Remark 1.1 and Theorem 1.3 \cite{DaOrVi}) that the stationary Airy sheet exists and is unique in law. The parabolic Airy sheet satisfies a version of the 1:2:3 scaling with respect to metric composition. For each $\gamma>0$ let $S_\gamma$ denote the diffusive scaling transform, which we will apply to real functions of one or two variables:
$$S_\gamma\ff(x):=\gamma^{-1}\ff(\gamma^2 x)\,\mbox{ and }\,S_\gamma\ff(x,y):=\gamma^{-1}\ff(\gamma^2 x,\gamma^2 y)\,.$$
Define the parabolic Airy sheet $\cL_s$ of scale $s>0$ by 
$$\cL_s(x,y):=S_{s^{-1}}\cL(x,y)=s\cL(x/s^2,y/s^2)\,.$$    
Then
\begin{equation*}\label{123MetComp}
\cL_r(x,y)\stackrel{dist.}{=}\max_{z\in\R}\left\{\cL^{(1)}_s(x,z)+\cL^{(2)}_t(z,y)\right\}\,,\,\mbox{ with }r^3=s^3+t^3\,,
\end{equation*}
(as random functions) where $\cL^{(1)}_s$ and $\cL^{(2)}_t$ are two independent copies of the parabolic Airy sheet of scales $s,t>0$, respectively. (For the parabolic Airy sheet \eqref{PaAiSh} we have a true maximum!)
\newline

To introduce the directed landscape we consider an oriented four-dimensional  parameter space defined as 
$$\R^4_{\uparrow}:=\left\{(x,s;y,t)\in\R^4\,:\,s<t\right\}\,.$$ 
Coordinates $s$ and $t$ represents time while coordintes $x$ and $y$ represents space. In the next we follow Definition 10.1 \cite{DaOrVi} to introduce the directed landscape. By Theorem 10.9 \cite{DaOrVi}, the directed landscape exists and is unique in law.

\begin{Definition}\label{DirLand}
The directed landscape is a random continuous function $\cL:\R_\uparrow^4\to\R$ that satisfies the following properties.
\begin{itemize}
\item Airy sheets marginals: for each $t\in\R$ and $s>0$ we have
\begin{equation}\label{DirSheet}
\cL(\cdot,t;\cdot,t+s)\stackrel{dist.}{=}\cL_s(\cdot,\cdot)\,.
\end{equation}
\item Independent increments: if $\{(t_i,t_i+s_i)\,:\,i=1,\dots,k\}$ is a  collection of disjont intervals then $\{\cL(\cdot,t_i;\cdot,t_i+s_i)\,:\,i=1,\dots,k\}$ is a collection of independent random functions. 
\item Metric composition: almost surely
\begin{equation}\label{MetComp}
\cL(x,r;y,t)=\max_{z\in\R}\left\{\cL(x,r;z,s)+\cL(z,s;y,t)\right\}\,,\,\forall\,(x,s;y,t)\in\R_{\uparrow}^4\mbox{ and }s\in(r,t)\,.
\end{equation}
\end{itemize}
\end{Definition}

Dauvergne, Ortmann and Virag \cite{DaOrVi} showed that the directed landscape describes the full space and time scaling limit of the fluctuations of the Brownian last-passage percolation model (Theorem 1.5 \cite{DaOrVi}). By setting $(x,s)_n:=(s+2x/n^{1/3},-\lfloor sn \rfloor)$, they proved that there exists a coupling between the directed landscape and the Brownian last-passage percolation model such that 
\begin{equation}\label{FullScal}
\B^{(n)}\left((x,s)_n\to(y,t)_n\right)=2(t-s)\sqrt{n}+2(y-x)n^{1/6}+ n^{-1/6}\left(\cL+o_n\right)(x,s;y,t)\,,
\end{equation}
where $\B^{(n)}(\dots,\B^{(n)}_{-1},\B^{(n)}_0,\B^{(n)}_1,\dots)$ is a sequence of Brownian motions and $o_n$ is a random function asymptotically small in the sense that for each compact $K\subseteq\R^4_\uparrow$ there exists $a>1$ such that $\E\left(a^{\sup_{K}o_n}\right)\to 1$ as $n\to\infty$. 
\newline

The directed landscape induces a semi-group evolution which takes into account the metric composition \eqref{MetComp}. The space $\sfUC$ is defined below. As our initial data, we incorporate (generalized) functions that might take value $-\infty$. 
\begin{Definition}\label{DefUC}
We say that a function $\ff:\R\to [-\infty,\infty)$ is upper semicontinuous if  
$$\limsup_{x\to y}\ff(x)\leq \ff(y)\,.$$ 
Let $\sfUC$ denote the space of upper semicontinuous generalized functions $\ff:\R\to [-\infty,\infty)$ with $\ff(x)\leq C_1|x|+C_2$ for all $x\in\R$, for some $C_1,C_2<\infty$, and $\ff(x)>-\infty$ for some $x\in\R$.
\end{Definition}

A canonical example of a (generalized) upper semicontinuous function that will be consider here several time is 
\begin{equation}\label{Dirac}
\fd_x(z)=\left\{\begin{array}{ll}0 & \mbox{ for } z=x\\-\infty &\mbox{ for } z\neq x\,.\end{array}\right. 
\end{equation}
The state space $\sfUC$ can be endowed with the topology of local convergence turning it into a Polish space (Section 3.1 \cite{MaQuRe}). 
\begin{prop}\label{KPZDef}
Let $\fh\in\sfUC$. Then a.s. for all $0<s<t$ and $x\in\R$ the random function $z\in\R\mapsto\fh(z)+\cL(z,s;x,t)$ attains it maximum and the process 
\begin{equation}\label{EvoDef}
\fh_{s,t}(x;\fh):=\max_{z\in\R}\left\{\fh(z)+\cL(z,s;x,t)\right\}\,,
\end{equation} 
defines a Markov semi-group acting on $\sfUC$, i.e. 
$$\fh_{r,t+s}(\cdot;\fh)=\fh_{t,t+s}(\cdot;\fh_{r,t})\,.$$
From now on we denote $\fh_t\equiv\fh_{0,t}$.
\end{prop}

Proposition \ref{KPZDef} follows from Proposition \ref{PropBasic}, which will be proved in the next section. Notice that, by independence of increments (Definition \ref{DirLand}), $\fh_{r,t}(\cdot;\fh)$ and $\cL(\cdot,t;\cdot,t+s)$ are independent. The directed landscape can be recovered in terms of semi-group $\fh_{s,t}$ by choosing a proper initial condition \eqref{Dirac}:
\begin{equation}\label{DirEvo}
\cL(x,s;y,t)=\fh_{s,t}(y;\fd_x)\,.
\end{equation}
The KPZ fixed point satisfies the so called 1-2-3 scaling invariance:
\begin{equation}\label{123}
S_{\gamma^{-1}}\fh_{\gamma^{-3}t}(\cdot;S_\gamma\fh)\stackrel{dist.}{=}\fh_{t}(\cdot;\fh)\,.
\end{equation}
Furthermore, if we set 
\begin{equation}\label{2BM_1}
\fb\equiv\mbox{ two-sided Brownian motion with diffusion coefficient $2$}\,,
\end{equation}
then 
\begin{equation}\label{stat_1}
\Delta \fh_t(\cdot;\fb^{\mu})\stackrel{dist.}{=}\fb^\mu(\cdot)\,,\,\mbox{ for all }\,t\geq 0\,,
\end{equation}
where $\Delta\ff(x):=\ff(x)-\ff(0)$ and $\fb^{\mu}(x):=\mu x+\fb(x)$.
\newline

The transition probabilities of the semi-group $\fh_{s,t}$ were computed by Matetski, Quastel and Remenik \cite{MaQuRe}, and we give a brief description as follows. Notice that the collection composed by cylindrical subsets of $\sfUC$,    
$$\mathrm{Cy}(\vec{x},\vec{a}):=\left\{\ff\in\sfUC\,:\,\ff(x_1)\leq a_1,\dots,\ff(x_m)\leq a_m\right\}\,\mbox{ for }\vec{x},\vec{a}\in\R^m\,,$$
is a generating sub-algebra for the Borel $\sigma$-algebra over $\sfUC$. The KPZ fixed point $\left(\fh_t(\cdot)\,,\,t\geq 0\right)$ is the unique time homogenous Markov process taking values in $\sfUC$ with transition probabilities given by the extension from the cylindrical sub-algebra to the Borel sets of 
\begin{equation}\label{DefKPZ}
\P\Big(\fh_t\in \mathrm{Cy}(\vec{x},\vec{a})\mid\fh_0=\fh\Big)=\det\left(\II-\KK^{\fh}_{t,\vec{x},\vec{a}}\right)_{\bL^2(\{x_1,\dots,x_m\}\times\R)}\,.
\end{equation}
On the right hand side of \eqref{DefKPZ} we have a Fredholm determinant of the operator $\KK^{\fh}_{t,\vec{x},\vec{a}}$, whose definition we address to \cite{MaQuRe} ($\II$ is the identity operator). From this formula one can recover several of the classical Airy processes by starting with special profiles for which the respective operators $K$ are explicit (see Section 4.4 of \cite{MaQuRe}). For instance, the Airy$_2$ process $\cA(\cdot)=\fh(\cdot;\fd)$ is defined by taking the initial profile $\fh=\fd$ where $\fd(0)=0$ and $\fd(z)=-\infty$ for all $z\neq 0$.
\newline

\subsection{Space H\"older Regularity and Brownian Behaviour}
Using kernel estimates for discrete approximations of the integral operator in \eqref{DefKPZ}, Matetski, Quastel and Remenik \cite{MaQuRe} proved that $\fh_t$ has H\"older $1/2-$ regularity in space (Theorem 4.13 \cite{MaQuRe}), and also that $S_{\sqrt{\epsilon}}\Delta\fh_t$ converges to $\fb$, as $\epsilon\to 0^+$, in terms of finite dimensional distributions (Theorem 4.14 \cite{MaQuRe}). Functional convergence was proved by Pimentel \cite{Pi} for several versions of Airy processes, which are obtained from the fundamental initial profiles $\fh\equiv\fd_0$, $\fh\equiv 0$ and $\fh\equiv \fb$, and stronger forms of local Brownian behaviour were proved by Corwin and Hammond \cite{CoHa} and Hammond \cite{Ha}. Here we use geometrical properties related to \eqref{EvoDef} to control space regularity of $\fh_{t}$. 
\newline

Let $\beta\in[0,1]$ and define the H\"older semi-norm of a real function $\ff:\R\to\R$ as  
$$\| \ff\|_{\beta,[-a,a]}:=\sup\left\{\frac{|\ff(x)-\ff(y)|}{|x-y|^\beta}\,:\,x,y\in[-a,a]\,\mbox{ and }x\neq y\right\}\,.$$
\begin{thm}\label{Holder_1}
Fix $a,t>0$ and $\beta\in[0,1/2)$. Then  
\begin{equation}\label{EqHolder}
\P\left(\| \fh_t\|_{\beta,[-a,a]}<\infty\right)=1\,.
\end{equation}
Furthermore,
\begin{equation}\label{EqBrown}
\lim_{\epsilon\to 0^+}S_{\sqrt{\epsilon}}\Delta\fh_t(\cdot)\stackrel{dist.}{=}\fb(\cdot)\,,
\end{equation}
where the distribution of $\fb$ is given by \eqref{2BM_1}.
\end{thm} 

For $\ff:\R^2\to\R$, define the H\"older semi-norm as follows
$$\| \ff\|_{\beta,[-a,a]^2}:=\sup\left\{\frac{|\ff(x_1,x_2)-\ff(y_1,y_2)|}{|(x_1,x_2)-(y_1,y_2)|_\infty^\beta}\,:\,(x_1,x_2),(y_1,y_2)\in[-a,a]\,\mbox{ and }(x_1,x_2)\neq (y_1,y_2)\right\}\,.$$
Denote  
$$\Delta\ff(x,y):=\ff(x,y)-\ff(0,0)\,,$$
and let $\cB(x,y):=\fb_1(x)+\fb_2(y)$, where $\fb_1$ and $\fb_2$ are two independent copies of \eqref{2BM_1}.   

\begin{thm}\label{AirySheet}
Consider the stationary Airy sheet and $\beta\in[0,1/2)$. Then 
\begin{equation*}\label{AiHolder}
\P\left(\|\cA\|_{\beta,[-a,a]^2}<\infty\right)=1\,.
\end{equation*}
Furthermore \footnote{Convergence in terms of a sequence of random elements in the space of continuous scalar fields on a fixed compact subset of $\R^2$, endowed with the uniform metric.}
$$\lim_{\epsilon\to 0^+}S_{\sqrt{\epsilon}}\Delta \cA(\cdot,\cdot) \stackrel{dist.}{=}\cB(\cdot,\cdot)\,.$$
In view of \eqref{PaAiSh} and \eqref{DirSheet}, we also have that 
$$\lim_{t\to \infty}\Delta\cL(\cdot,0;\cdot,t)\stackrel{dist.}{=}\cB(\cdot,\cdot)\,.$$
\end{thm}

\subsection{Brownian Long Time Behaviour}
From \eqref{123}, one can see that the long time behaviour of $\Delta\fh_t$ can be written in terms of the local space behaviour of $\Delta\fh_1$ (take $\gamma=t^{1/3}$), which allows one to obtain long time convergence (in terms of finite dimensional distributions) from the local convergence to Brownian motion, as soon as $S_\gamma\fh$ converges in distribution in $\sfUC$ as $\gamma\to\infty$ (Theorem 4.15 \cite{MaQuRe}). Based on the same geometrical tools to study the space regularity of the KPZ fixed point, we will prove long time convergence of the KPZ fixed.
\begin{Assump}\label{Assump}
There exist $c>0$ and a real function $\psi$ such that for all $\gamma\geq c$ and $r\geq 1$ 
\begin{equation}\label{Assump1}
\P\left(\,S_\gamma\fh(z)\leq r |z|\,,\,\forall\,|z|\geq 1\,\right)\geq 1-\psi(r)\,\,\mbox{ and }\,\,\lim_{r\to\infty}\psi(r)=0\,.
\end{equation}  
\end{Assump}

\begin{thm}\label{Coupling}
Let $a,t,\eta>0$ and set $r_t:=\sqrt[4]{t^{2/3}a^{-1}}$. Under Assumption \ref{Assump}, where $\fb$ \eqref{2BM_1} and $\fh$ are sample independently, there exists a real function $\phi$, which does not depend on $a,t,\eta>0$, such that for all $t\geq \max\{c^3,a^{3/2}\}$ and $\eta>0$ we have
\begin{equation}\label{Coupling1}
\P\left(\sup_{x\in[-a,a]}|\Delta\fh_t(x;\fh)-\Delta\fh_t(x;\fb)|>\eta \sqrt{a}\right)\leq \phi\left(r_t\right)+\frac{1}{\eta r_t}\,\,\mbox{ and }\,\,\lim_{r\to \infty}\phi(r)=0\,.
\end{equation}
In particular, if $\lim_{t\to\infty}a_{t} t^{-2/3}=0$ then 
$$\lim_{t\to\infty}\P\left(\sup_{x\in[-a_t,a_t]}|\Delta\fh_t(x;\fh)-\Delta\fh_t(x;\fb)|>\eta \sqrt{a_t}\right)=0\,.$$
Since $S_{\sqrt{a_t}}\Delta\fh_t(\cdot;\fb)\stackrel{dist.}{=}\fb(\cdot)$, we also have that 
$$\lim_{t\to\infty}S_{\sqrt{a_t}}\Delta\fh_t(\cdot;\fh)\stackrel{dist.}{=}\fb(\cdot)\,.$$
\end{thm}

\begin{rem}\label{Decay}
For deterministic $\fh(x)=x^{\zeta}$, for $\zeta\in[0,1]$, we have that $S_\gamma\fh(x)=\gamma^{2\zeta-1}x^\zeta$. If $\zeta\in[0,1/2]$, then $\fh$ does satisfy \eqref{Assump1}, while for $\zeta\in(1/2,1]$ it does not. We use assumption \eqref{Assump1} to ensure that, for all large values of $t$,
\begin{equation}\label{Assump2}
\P\left(|Z_t(\pm a;\fh)|>rt^{2/3}\right)\leq \phi_1(r) \to 0\,, \mbox{ as }r\to\infty\,,
\end{equation}
where $Z_t(x;\fh)$ is the rightmost $z\in\R$ to attain the maximum \eqref{EvoDef}, and $\phi_1$ is a real function that does not depend on $a>0$ or $t>0$ (Lemma \ref{KPZLocalization2}). If one can prove \eqref{Assump2}, based on possible different assumptions, then \eqref{Coupling1} will follow as well.  
\end{rem}

\begin{rem}\label{ErgDrift}
Theorem \ref{Coupling} does not imply immediately that the only spatially ergodic (in terms of its increments) and time invariant process with zero drift is $\fb$. This would follow as soon as one can verify \eqref{Assump1} or \eqref{Assump2} for such a process.
\end{rem}

\section{Geometry, Comparison and Attractiveness}\label{Directed}

Given an upper semicontinuous function $\ff$ such that   
\begin{equation}\label{decay}
\lim_{|z|\to\infty}\ff(z)=-\infty\,,
\end{equation}
then the supremum of $\ff(z)$ over $z\in\R$ is indeed a maximum, i.e. $\exists\,Z\in\R$ such that $\ff(Z)\geq \ff(z)$ for all $z\in\R$. Additionally, the set  
$$\argmax_{z\in\R}\ff(z):=\left\{Z\in\R\,:\,\ff(Z)=\max_{z\in\R}\ff(z)\right\}\,.$$ 
is compact. Since with probability one, for all $\fh\in \sfUC$, $\fh(z)+\cL(z,s;x,t)$ satisfies \eqref{decay}, for all $s<t$ and $x\in\R$, (due to the parabolic drift \eqref{PaAiSh}) we can use these aforementioned facts to study the semi-group evolution \eqref{EvoDef}. 
\newline

We call a continuous path $\cP:[r,t]\to\R$ a geodesic between the space-time points $(x,r)$ and $(y,t)$ if $\cP(r)=x$, $\cP(t)=y$ and for $s\in(r,t)$
\begin{equation}\label{DirGeo}
\cL(x,r;y,t)=\cL\left(x,r;\cP(s),s\right)+\cL\left(\cP(s),s;y,t\right)\,
\end{equation}
Define $\cP_{x,r}^{y,t}(r)=x$, $\cP_{x,r}^{y,t}(t)=y$ and 
$$\cP_{x,r}^{y,t}(s):=\max\argmax_{z\in\R}\left\{\cL(x,r;z,s)+\cL(z,s;y,t)\right\}\mbox{ for }s\in(r,t)\,.$$
By Lemma 13.3 \cite{DaOrVi}, almost surely, $\cP_{x,r}^{y,t}$ is a geodesic for every $(x,r)$ and $(y,t)$. We also identify the geodesic path (or function) $\cP$ with its graph $\{(\cP(s),s)\,:\,s\in[r,t]\}$ in order to handle intersection points between different paths. For each $\fh\in\sfUC$, $0<t$ and $x\in\R$, let
\begin{equation}\label{ArgDef}
Z_{t}(x;\fh):=\max\argmax_{z\in\R}\left\{\fh(z)+\cL(z,s;x,t)\right\}\,.
\end{equation}
\begin{prop}\label{PropBasic}
Almost surely $\fh_{t}$ and $Z_{t}$ are a well defined real functions for which we have the following properties.
\begin{itemize}
\item[(i)] $\fh_t(x)=\fh\left(Z_t(x)\right)+\cL\left(Z_t(x),0;x,t\right)$. 
\item[(ii)] For every $w\in\R$ and $u\in[0,t)$,
$$\fh_t(x)\geq \fh_u(w)+\cL(w,u;x,t)\,.$$ 
\item[(iii)] For every $(w,u)\in\cP_{Z_t(x),0}^{x,t}$,  
$$\fh_t(x)=\fh_u(w)+\cL(w,u;x,t)\mbox{ and }\fh_u(w)=\fh(Z_t(x))+\cL(0,Z_t(x);w,u)\,.$$ 
\item[(iv)] For fixed $t>0$, $Z_t(x)$ is a nondecreasing function of $x\in\R$.
\item[(v)] $\fh_{s,t}$ defines a semi-group: $\fh_{t+s}(\cdot;\fh)=\fh_{t,t+s}(\cdot;\fh_t)$, i.e.
$$\fh_{t+s}(x;\fh)=\max_{z\in\R}\left\{\fh_t(z;\fh)+\cL(z,t;x,t+s)\right\}\,,\,\forall\,x\in\R\,.$$
\end{itemize}

\end{prop}

\noindent\paragraph{\bf Proof} By compactness, $Z_t(x)\in\argmax_{z\in\R}\left\{\fh(z)+\cL(z,0;x,t)\right\}$, which implies (i). Now we use \eqref{EvoDef} and \eqref{MetComp} to get (ii): for any $z,w\in\R$ and $u\in (0,t)$,
$$\fh_t(x)\geq\fh(z)+\cL(0,z;x,t)\geq \fh(z)+\cL(0,z;u,w)+\cL(w,u;x,t)\,,$$
and hence
$$\fh_t(x)\geq \fh_u(w)+\cL(w,u;x,t)\,.$$
By (i) and \eqref{DirGeo}, if $w=\cP_{Z_t(x),0}^{x,t}(u)$ then 
\begin{eqnarray*}
\fh_t(x)&=&\fh(Z_t(x))+\cL(Z_t(x),0;x,t)\\
&=&\fh(Z_t(x))+\cL(Z_t(x),0;w,u)+\cL(w,u;x,t)\\
&\leq & \fh_u(w)+\cL(w,u;x,t)\,,
\end{eqnarray*}
and thus, by (ii), $\fh_t(x)= \fh_u(w)+\cL(w,u;x,t)$ and $\fh_u(w)=\fh(Z_t(x))+\cL(0,Z_t(x);w,u)$, which concludes the proof of (iii). To prove (iv), assume that $Z_t(y)<Z_t(x)$ for some $x<y$. Then $\cP_{Z_t(y),0}^{y,t}$ and $\cP_{Z_t(x),0}^{x,t}$ intersects at some space-time point $(w,u)$. By (iii), we have that
$$\fh_t(y)=\fh_u(w)+\cL(w,u;y,t)\,\,\mbox{ and }\,\,\fh_u(w)=\fh(Z_t(x))+ \cL(Z_t(x),0;w,u)\,.$$
This shows that
\begin{eqnarray*}
\fh_t(y)&=&\fh_u(w)+\cL(w,u;y,t)\\
&=&\fh(Z_t(x))+ \cL(Z_t(x),0;w,u)+\cL(w,u;y,t)\\ 
&\leq &\fh(Z_t(x))+ \cL(Z_t(x),0;y,t)\,,
\end{eqnarray*}
where we use the metric composition \eqref{MetComp} for the last inequality. Hence, $Z_t(x)$ is also a location that attains the maximum for $\fh_t(y)$,   
which leads to a contradiction since we assumed that $Z_t(y)<Z_t(x)$ and $Z_t(y)$ is the rightmost point to attain the maximum. The the semi-group property (v) follows directly  item (iii). 

\hfill$\Box$\\

\begin{prop}[Argmax Comparison]\label{ArgmaxComparison}
If $x<y$ and $Z_t(y; \fh)\leq Z_t(x; \tilde\fh)$ then 
$$\fh_t(y;\fh)-\fh_t(x;\fh)\leq \fh_t(y;\tilde\fh)-\fh_t(x;\tilde\fh)\,.$$
\end{prop}

\noindent\paragraph{\bf Proof } Denote $z\equiv Z_t(y;\fh)$ and $\tilde z\equiv Z_t(x;\tilde\fh)$. By assumption, $x<y$ and $z\leq \tilde z$, and hence there exists $(w,u)\in\cP_{z,0}^{y,t}\cap\cP_{\tilde z,0}^{x,t}$. Since $(w,u)\in\cP_{\tilde z,0}^{x,t}$, by (iii)-Proposition \ref{PropBasic},
$$\fh_t(x;\tilde\fh)= \fh_u(w;\tilde \fh)+\cL(w,u;x,t)\,,$$
and, by (ii)-Proposition \ref{PropBasic},
$$\fh_t(y;\tilde\fh)  \geq  \fh_u(w;\tilde\fh)+\cL(w,u;y,t) \,,$$
that yields to 
$$\fh_t(y;\tilde\fh) - \fh_t(x;\tilde\fh)\geq \cL(w,u;y,t) - \cL(w,u;x,t)\,.$$
Now $(w,u)\in\cP_{z,0}^{y,t}$ and by using Proposition \ref{PropBasic} as before, we have  
$$\fh_t(y;\fh)=\fh_u(w;\fh)+\cL(w,u;y,t)\,\mbox{ and }\,\fh_t(x;\fh)\geq \fh_u(w;\fh)+\cL(w,u;x,t)\,,$$
which implies that 
$$\fh_t(y;\fh)-\fh_t(x;\fh)\leq \cL(w,u;y,t) - \cL(w,u;x,t)\,,$$
and therefore 
$$\fh_t(y;\fh)-\fh_t(x;\fh)\leq\fh_t(y;\tilde\fh) - \fh_t(x;\tilde\fh)\,.$$

\hfill$\Box$\\ 

\begin{prop}[Attractiveness]\label{Attractiveness}
If $\fh(y)-\fh(x)\leq\tilde\fh(y)-\tilde\fh(x)$ for all $x<y$ then
$$\fh_t(y;\fh)-\fh_t(x;\fh)\leq \fh_t(y;\tilde\fh)-\fh_t(x;\tilde\fh)\,\,\,\forall\,x<y\,,\forall\,t\geq 0\,.$$ 
\end{prop}

\noindent\paragraph{\bf Proof} Denote again $z\equiv Z_t(y;\fh)$ and $\tilde z\equiv Z_t(x;\tilde\fh)$. If $z \leq \tilde z$ then 
$$\fh_t(y;\fh)-\fh_t(x;\fh)\leq \fh_t(y;\tilde\fh)-\fh_t(x;\tilde\fh)\,,$$
by Proposition \ref{ArgmaxComparison}. If  $z>\tilde z$ then, by (i)-Proposition \ref{PropBasic},
\begin{eqnarray*}
\fh_t(y;\tilde\fh) - \fh_t(x;\tilde\fh)  &-&\Big(\fh_t(y;\fh) - \fh_t(x;\fh) \Big)\\
&=&\fh_t(y;\tilde\fh)-\big(\tilde\fh(\tilde z)+\cL(\tilde z,0;x,t)\big)-\Big(\big(\fh(z)+\cL(z,0;y,t)\big)-\fh_t(x;\fh)\Big)\\
&=&\fh_t(y;\tilde\fh)-\big(\tilde\fh(z)+\cL(z,0;y,t)\big)+\Big(\fh_t(x;\fh)- \big(\fh(\tilde z)+\cL(\tilde z,0;x,t)\big)\Big)\\
&+&\big(\tilde\fh(z)-\tilde\fh(\tilde z)\big)-\big(\fh(z)-\fh(\tilde z)\big) \,.
\end{eqnarray*}
Thus, by \eqref{EvoDef},
$$\fh_t(y;\tilde\fh)-\big(\tilde\fh(z)+\cL(z,0;y,t)\big)\geq 0\,\,\mbox{ and }\,\,\fh_t(x;\fh)- \big(\fh(\tilde z)+\cL(\tilde z,0;x,t)\big)\geq 0\,,$$
while, by assumption, 
$$\big(\tilde\fh(z)-\tilde\fh(\tilde z)\big)-\big(\fh(z)-\fh(\tilde z)\big)\geq 0\,,$$
since $z>\tilde z$. 

\hfill$\Box$\\

\subsection{Uniqueness of the Argmax}
We finish this section by pointing out how the ideas in \cite{Pi1} can be combined with the fact that the Airy$_2$ process is locally absolutely continuous with respect to Brownian motion \cite{CoHa}, to prove a.s. uniqueness of the location of the maxima in \eqref{EvoDef}. Since $\fh(z)+\cL(z,s;x,t)$ satisfies \eqref{decay}, it is enough to prove uniqueness of the location of the maximum restrict to a compact set. On the other hand,  $\{\cL(z,s;x,t)\,:\,z\in\R\}$ is distributed as a rescaled Airy$_2$ process minus a parabola (for fixed $x\in\R$ and $0<s<t$), which is locally absolutely continuous with respect to Brownian motion \cite{CoHa}. Therefore, uniqueness of the location of the maxima in \eqref{EvoDef} follows from the next proposition, which is similar to Theorem 2 \cite{Pi1}, combined with Lemma 2 \cite{Pi1}\footnote{Lemma 2 \cite{Pi1} shows that $m(a)$ is differentiable at $a=0$ if $\ff$ is a sum of a deterministic function $\fh$ with a Brownian motion.}.   

\begin{prop}\label{UniqArgMax}
Let $K\subseteq \R$ be a compact set and $\ff:K\to\R$ be a random upper semicontinuous function. Denote $\ff^a(z):=\ff(z)+az$, $M(\ff):=\max_{z\in K}\ff(z)$ and 
$$m(a)=\E\left(M\left(\ff^a\right)\right)\,.$$
 Then a.s. there exists a unique $Z\in K$ such that $M(\ff)=\ff(Z)$  
if and only if $m(a)$ is differentiable at $a=0$. Furthermore, in this case,
$$m'(0)=\E Z\,.$$  
\end{prop}

\noindent\paragraph{\bf Proof} The first part of the proof is merely analytic and we follow the proof of Lemma 1 \cite{Pi1}, where $\ff$ was assumed to be continuous. There are two fundamental steps where we used continuity that needs to be adapted to upper semicontinuous functions. Denote  
$$Z_1(\ff):=\inf\argmax_{z\in K} \ff(z)\,\,\mbox{ and }\,\,Z_2(\ff):=\sup\argmax_{z\in K} \ff(z)\,.$$ 
For simple notation we put $M^a\equiv M(\ff^a)$, $M\equiv M(\ff)$, $Z^a_i\equiv Z_i(\ff^a)$ and finally $Z_i\equiv Z_i(\ff)$. The first step in \cite{Pi1} was to argue that $M=\ff(Z_i)$ and $M^a=\ff(Z_i^a)+a Z_i^a$. But for a upper semicontinuous function, $\argmax_{z\in K} \ff(z)$ is a compact set, and then $Z_1(\ff),Z_2(\ff)\in \argmax_{K} \ff(z)$ (which also holds for $\ff^a$). Thus, we can conclude that 
\begin{equation}\label{step1}
M+aZ_i=\ff(Z_i)+aZ_i\leq M^a=\ff(Z_i^a)+a Z_i^a \leq M+ a Z_i^a\,.
\end{equation}
The second step in \cite{Pi1} was to prove that    
\begin{equation}\label{step2}
\lim_{a\to 0^-}Z_1^a=Z_1\,\mbox{ and }\,\lim_{a\to 0^+}Z_2^a=Z_2\,.
\end{equation}
Indeed, by \eqref{step1}, we have that $Z_1^a\leq Z_1$ for all $a<0$, and if the convergence of $Z_1^a$ to $Z_1$ does not hold then, by compactness of $K$, we can find $\tilde Z_1\in K$, $\delta>0$  and a sequence $a_n\to 0^-$ such that $\lim_{n\to\infty}Z_1^{a_n}=\tilde Z_1$ and $\tilde Z_1\leq Z_1-\delta$. But by \eqref{step1}, we also have that 
$$0\leq a\left(Z_i^a-Z_i\right)-\left(\ff(Z_i)-\ff(Z_i^a)\right)\,,\mbox{ for $i=1,2$}\,,$$
and thus (first inequality)
$$\ff(Z_1)\leq \limsup_{n}\ff(Z_1^{a_n})\leq \ff(\tilde Z_1)\,,$$
where we use upper semicontinuity in the second inequality. But this is a contradiction, since $Z_1$ is the leftmost location to attain the maximum, and hence $\lim_{a\to 0^-}Z^a_1=Z_1$.  Since $Z_2^a\geq Z_2$ for all $a>0$, the proof of $\lim_{a\to 0^+}Z^a_2=Z_2$ is analogous. By \eqref{step1} again,
$$0\leq (M^a-M)-aZ_i\leq a(Z_i^a-Z_i)\,,$$
which implies that 
$$0\geq \frac{M^a-M}{a} - Z_1\geq Z^a_1-Z_1\geq -diam(K)\,,\,\mbox{ for }a<0\,,$$
and 
$$ 0\leq \frac{M^a-M}{a} - Z_2\leq Z^a_2-Z_2\leq diam(K)\,,\,\mbox{ for }a>0\,,$$
where $diam(K)$ denotes the diameter of $K$. Since the location of the maxmimum is a.s. unique if and only if $\E\left(Z_1\right)=\E\left(Z_2\right)$ (now we have a random $\ff$), using the inequalities above, \eqref{step2} and dominated convergence, we see that the location of the maximum of $\ff$ is a.s. unique if and only if $m(a)$ is differentiable at $a=0$:
$$\E\left(Z_1\right)=\E\left(Z_2\right)\,\Longleftrightarrow\,\lim_{a\to 0^-}\frac{m(a)-m(0)}{a}=\lim_{a\to 0^+}\frac{m(a)-m(0)}{a}\,.$$

\hfill$\Box$\\ 

\section{Proof of the Theorems}
\subsection{KPZ Localization of the Argmax}
A key step to use comparison (Proposition \ref{ArgmaxComparison}) relies on the control of  $Z_t(x;\fh)$ (recall \eqref{ArgDef}) as a function of $\fh$, $x$ and $t$. Let $X$ be the closest point to the origin such that $\fh(X)>-\infty$ (if a tiebreak occurs we pick the nonnegative one). By assumption, $X\in\R$ is a well defined random variable. Since the location of a maximum is invariant under vertical shifts of $\fh$, if we want to control the location of the maximum, we can assume without loss of generality that $\fh(X)=0$. By the symmetries (i)-(ii)-(iii), for fixed values of $x\in\R$ and $t>0$,
\begin{equation}\label{ArgDist}
Z_t(x;\fh)\stackrel{dist.}{=} t^{2/3}Z_1(0;S_{\gamma_t}T_x\fh)+x\,,
\end{equation}
where $\gamma_t:=t^{1/3}$. By \eqref{ArgDist}, for $x\in[-a,a]$,
\begin{equation}\label{ArgTailCont1}
\P\left(|Z_t(x;\fh)|>rt^{2/3}\right)\leq\P\left(|Z_1(0;S_{\gamma_t}T_x\fh )|>r-|x|t^{-2/3}\right)\leq \P\left(|Z_1(0;S_{\gamma_t}T_x\fh )|>r-at^{-2/3}\right)\,.
\end{equation}
The right hand side of \eqref{ArgTailCont1} is bounded by
$$\P\left(\max_{|z|>r-at^{-2/3}}\left\{S_{\gamma_t}T_x\fh(z) +\cA(z)-z^2\right\}=\max_{z\in\R}\left\{S_{\gamma_t}T_x\fh (z)+\cA(z)-z^2\right\}\right)\,,$$
where $\cA(z):=\cA(z,0)$. If we take $z=\gamma_t^{-2}(X-x)$ we get that $S_{\gamma_t}T_x\fh (z)=\fh(X)=0$, and  the right hand side of \eqref{ArgTailCont1} is bounded by 
\begin{equation}\label{ArgTailCont2}
\P\left(\max_{|z|>r-at^{-2/3}}\left\{S_{\gamma_t}T_x\fh (z)+\cA(z)-z^2\right\}\geq \cA\left(\gamma_t^{-2}(X-x)\right)-\left(\gamma_t^{-2}(X-x)\right)^2\right)\,.
\end{equation}
In the next lemmas we will use that the Airy$_2$ process $\{\cA(z)\,:\,z\in\R\}$ is stationary and  independent of $X$, which implies that $\cA\left(\gamma_t^{-2}(X-x)\right)\stackrel{dist.}{=}\cA(0)$, and  we can split the probability in \eqref{ArgTailCont2} as    
\begin{equation}\label{ArgTailCont3}
\P\left(\max_{|z|>r-at^{-2/3}}\left\{S_{\gamma_t}T_x\fh (z)+\cA(z)-z^2\right\}\geq -L \right)+\P\Big(\cA(0)-\gamma_t^{-4}(X-x)^2\leq -L\Big)\,,
\end{equation}
for any choice of $L>0$. 

\begin{lem}\label{KPZLocalization1}
Let $a,t>$ be fixed. For every $\fh\in\sfUC$
$$\lim_{r\to\infty}\P\left(\,|Z_t(\pm a;\fh)|>r\,\right)=0\,.$$
\end{lem}

\noindent\paragraph{\bf Proof} For the sake of simplicity, we are going to prove it for $t=1$ and $a=1$. Let us pick $L_r=(r-1)^2/4$. Then 
$$\lim_{r\to\infty}\P\left(\cA(0)-(X-1)^2\leq -L_r\right)=0\,,$$
since the random variable $\cA(0)-(X-1)^2$ does not depend on $r$. By \eqref{ArgTailCont1}, \eqref{ArgTailCont2} and \eqref{ArgTailCont3}, we still need to prove that  
$$\lim_{r\to\infty}\P\left(\max_{|z|>r-1}\left\{T_1\fh (z)+\cA(z)-z^2\right\}\geq -L_r \right)=0\,.$$
If $r>2$ and $|z|>r-1$ then $|z|>r/2>1$ and $\frac{r}{4}|z|-z^2\leq -z^{2}/2$. Hence, if $T_1\fh (z)\leq \frac{r}{4}|z|$ then 
$$ T_1\fh (z)-z^2\leq\frac{r}{4}|z|-z^2\leq -z^{2}/2 \,,$$
which shows that 
$$\P\left(\max_{|z|>r-1}\left\{T_1\fh (z)+\cA(z)-z^2\right\}\geq -L_r\right)\leq \psi\left(\,r/4;T_1\fh\,\right)+\P\left(\max_{|z|>r-1}\left\{\cA(z)-\frac{z^2}{2}\right\}\geq -L_r \right)\,,$$
where 
$$\psi(r;\fh):=1-\P\left(\,\fh(z)\leq r |z|\,,\,\forall\,|z|\geq 1\,\right)\,.$$
By (b)-Proposition 2.13 \cite{CoLiWa}, there exist constants $c_1,c_2>0$ such that for all $r>c_1$, 
$$\P\left(\max_{|z|>r}\left\{A(z)-\frac{z^2}{2}\right\}>-\frac{r^2}{4}\right)\leq e^{-c_2 r^3}\,,$$
which concludes the proof of,
$$\lim_{r\to\infty}\P\left(\,|Z_1(1;\fh)|>r\,\right)=0\,,$$
as soon as we prove that,
$$\lim_{r\to\infty}\psi(r;\fh)=0\,\mbox{ for all }\fh\in\sfUC\,.$$
But for every probability measure on $\sfUC$, we have that 
$$\P\left(\exists\,r>0\,\mbox{ s. t. }\fh(z)\leq r(1+|z|)\,\,\forall\,z\in\R\right)=1\,,$$
and if $r_1<r_2$ then
$$\{\fh(z)\leq r_1(1+|z|)\,\,\forall\,z\in\R\}\subseteq\{\fh(z)\leq r_2(1+|z|)\,\,\forall\,z\in\R\}\,,$$ 
which implies that
$$\lim_{r\to\infty}\P\left(\fh(z)\leq r(1+|z|)\,\,\forall\,z\in\R\right)=1\,.$$
Since $\frac{r}{2}(1+|z|)\leq r|z|$ for all $|z|\geq 1$ we have that 
$$\{\fh(z)\leq \frac{r}{2}(1+|z|)\,\,\forall\,z\in\R\}\subseteq\{\fh(z)\leq r|z|\,\,\forall\,|z|\geq 1\}\,,$$ 
and therefore, $\lim_{r\to\infty}\psi(r;\fh)=0$. 

\hfill$\Box$\\

\begin{lem}\label{KPZLocalization2}
Under \eqref{Assump1}, there exists a real function $\phi_1$, which does not depend on $a>0$ or $t>0$, such that for all $t\geq \max\{c^3,a^{3/2}\}$ we have
$$\P\left(|Z_t(\pm a;\fh)|>rt^{2/3} \right)\leq \phi_1(r)\,\,\mbox{ and }\,\,\lim_{r\to \infty}\phi_1(r) = 0\,.$$
\end{lem}

\noindent\paragraph{\bf Proof} Pick again $L_r=(r-1)^2/4$ and $t\geq \max\{c^3,a^{3/2}\}$. Then (recall that $\gamma_t=t^{1/3}$)
$$\gamma_t^{-4}(X-a)^2\leq 2\gamma_t^{-4}\left(X^2+a^2\right)\leq 2\left(\frac{X^2}{c^{4}}+1\right)\,,$$
and thus,
$$ \P\left(\cA(0)\leq -L_r+\gamma_t^{-4}(X-a)^2\right)\leq \P\left(\cA(0)\leq -L_r+ 2\left(\frac{X^2}{c^4}+1\right)\right)\,.$$
The right hand side of the above inequality is a function of $r$ that does not depend on $a>0$ or $t>0$, and goes to zero as $r$ goes to infinity. To control the other term in \eqref{ArgTailCont3} we note that, if $S_{\gamma_t}\fh(z)\leq \frac{r}{4} |z|$ for all $|z|\geq 1$, then 
$$S_{\gamma_t}T_a\fh(z)=t^{-1/3}T_a\fh(t^{2/3}z)=S_{\gamma_t}\fh\left(z+at^{-2/3}\right)\leq \frac{r}{4}|z+at^{-2/3}|\leq \frac{r}{4}\left(|z|+at^{-2/3}\right)\,,$$ 
as soon as $|z+at^{-2/3}|\geq 1$. This needs to hold for all $|z|>r-at^{-2/3}$ in order to upper bound the maximum over all such $z$'s. However, for $r>3$ and $|z|>r-at^{-2/3}$ (recall that $t\geq \max\{c^3,a^{3/2}\}$) we certainly have that $|z+at^{-2/3}|\geq 1$. Therefore, if $S_{\gamma_t}\fh(z)\leq \frac{r}{4} |z|$ for all $|z|\geq 1$, then  
\begin{eqnarray*}
\max_{|z|>r-at^{-2/3}}\left\{S_{\gamma_t}T_a\fh (z)+\cA(z)-z^2\right\}&\leq & \max_{|z|>r-at^{-2/3}}\left\{ \frac{r}{4}\left(|z|+at^{-2/3}\right)+\cA(z)-z^2\right\}\\
&\leq &  \max_{|z|>r-1}\left\{ \frac{r}{4}\left(|z|+1\right)+\cA(z)-z^2\right\}\,.
\end{eqnarray*}
To ensure that $\frac{r}{4}\left(|z|+1\right)-z^2\leq-\frac{z^2}{2}$ for $|z|>r-1$ we take $r>4$. Thus, we can conclude that for $r>4$ and $t\geq \max\{c^3,a^{3/2}\}$ we have that 
\begin{eqnarray*}
\P\left(\max_{|z|>r-at^{-2/3}}\left\{S_{\gamma_t}T_a\fh (z)+\cA(z)-z^2\right\}\geq -L_r \right)&\leq & \psi(r/4)\\
&+&\P\left(\max_{|z|>r-1}\left\{\cA(z)-\frac{z^2}{2}\right\}\geq -L_r \right)\,.
\end{eqnarray*}

\hfill$\Box$\\ 

For $\mu\geq 0$ denote
\begin{equation}\label{DriftBM}
\fh^{\pm\mu}_t(\cdot)\equiv\fh_t(\cdot;\fb^{\pm\mu})\,,\mbox{ where }\fb^{\pm\mu}(z)=\pm\mu z+\fb(z)\,,
\end{equation}
and $\fb$ is given by \eqref{2BM_1}. Hence, for all $x<y$,
\begin{equation}\label{DriftBM_1}
\fb^{-\mu}(y)-\fb^{-\mu}(x)\leq \fb(y)-\fb(x) \leq \fb^{\mu}(y)-\fb^{\mu}(x)\,.
\end{equation}
Recall \eqref{ArgDef} and let 
$$Z^{\pm\mu}_t(x)=Z_t\left(x;\fb^{\pm\mu}\right)\,.$$
Then 
\begin{equation}\label{BMArgMax}
Z^{\mu}_t(x)\stackrel{dist.}{=} Z^{0}_t(x)+\frac{\mu}{2}t\stackrel{dist.}{=}t^{2/3}Z^{0}_1(0)+x+ \frac{\mu}{2}t=t^{2/3}\left(Z^{0}_1(0)+xt^{-2/3}+ \frac{\mu}{2}t^{1/3}\right)\,.
\end{equation}
The next step is to construct an event $E_t(\mu)$ where we can sandwich the local increments of $\fh_t$ in between the local increments $\fh_t^{\pm\mu}$, and this is the point where we use Proposition \ref{ArgmaxComparison}. Define the event 
\begin{equation}\label{CouplingEvent}
E_t(\mu)=\left\{Z_t(a;\fh)\leq Z_t^{+\mu}(-a)\,\mbox{ and }Z_t(-a;\fh)\geq Z_t^{-\mu}(a)\,\right\}\,.
\end{equation}
By (iv)-Proposition \ref{PropBasic}, on the event $E_t(\mu)$, for $x<y$ and $x,y\in[-a,a]$,
$$Z_t(y;\fh)\leq Z_t(a;\fh)\leq Z_t^{\mu}(-a)\leq Z_t^{\mu}(x)\,,$$
and 
$$Z_t^{-\mu}(y)\leq Z_t^{-\mu}(a)\leq Z_t(-a;\fh)\leq Z_t(x;\fh)\,.$$  
Therefore, by Proposition \ref{ArgmaxComparison}, on the event $E_t(\mu)$,
if $x<y$ and $x,y\in[-a,a]$, then 
\begin{equation}\label{BasicComparison}
\fh_t^{-\mu}(y)-\fh_t^{-\mu}(x)\leq \fh_t(y;\fh)-\fh_t(x;\fh)\leq \fh_t^{\mu}(y)-\fh_t^{\mu}(x)\,.
\end{equation}

\subsection{Proof of Theorem \ref{Holder_1}} 
We want to control the H\"older semi-norm for $\beta\in[0,1/2)$ (we omit the dependence on the domain and on the initial profile $\fh$),
$$\| \fh_t\|_{\beta}\equiv\| \fh_t\|_{\beta,[-a,a]}:=\sup_{x, y\in[-a,a]\,,\,x\neq y}\frac{|\fh_t(x)-\fh_t(y)|}{|x-y|^\beta}\,.$$
By \eqref{BasicComparison}, on the event $E_t(\mu)$ \eqref{CouplingEvent},
$$\| \fh_t\|_{\beta} \leq \max\left\{\| \fh^{\mu}_t\|_{\beta}\,,\,\| \fh^{-\mu}_t\|_{\beta}\right\}\,,$$
and hence,
\begin{eqnarray*}
\P\left(\| \fh_t\|_{\beta}>A\right)& \leq &\P\left(\left\{\| \fh_t\|_{\beta}>A\right\}\cap E_t(\mu)\right)+\P\left(E_t(\mu)^c\right)\\
&\leq &\P\left(\| \fh^{\mu}_t\|_{\beta}>A\right)+\P\left(\| \fh^{-\mu}_t\|_{\beta}>A\right)+\P\left(E_t(\mu)^c\right)\,.
\end{eqnarray*}
Since $\|\fh^{\pm\mu}_t\|_{\beta}= \|\Delta\fh^{\pm\mu}_t\|_{\beta}$ and $\Delta\fh^{\pm\mu}_t$ are drifted Brownian motions \eqref{stat_1}, 
$$ \lim_{A\to\infty}\P\left(\| \fh^{\pm\mu}_t\|_{\beta}>A\right)=\lim_{A\to\infty}\P\left(\| \Delta\fh^{\pm\mu}_t\|_{\beta}>A\right)=0\,,$$
which yields to 
$$0\leq \limsup_{A\to\infty} \P\left(\| \fh_t\|_{\beta}>A\right)\leq \P\left(E_t(\mu)^c\right)\,.$$
We picked $\mu>0$ arbitrary and
\begin{equation}\label{ControlComparison}
\P\left(E_t(\mu)^c\right)\leq\P\left(|Z_t(a;\fh)|>\frac{\mu}{4}t\right)+ \P\left(Z^{\mu}_t(-a)\leq \frac{\mu}{4}t\right)+\P\left(Z^{-\mu}_t(a)\geq -\frac{\mu}{4}t\right)\to 0\,,
\end{equation}
as $\mu\to\infty$, by \eqref{BMArgMax} and Lemma \ref{KPZLocalization1}. Therefore
$$\lim_{A\to\infty} \P\left(\| \fh_t\|_{\beta}>A\right)=0\,,$$
which implies \eqref{EqHolder}.

To prove convergence of 
$$S_{\sqrt{\epsilon}}\Delta\fh_t(x;\fh)=\epsilon^{-1/2}\left(\fh_t(\epsilon x)-\fh_t(0)\right)\,,$$
to Brownian motion \eqref{EqBrown}, we consider the event $E_t(\mu)$ \eqref{CouplingEvent} again with $a=1$ (we will choose $\mu$ later as a suitable function of $\epsilon$). Given a compact set $K\subseteq \R$ we take a $\epsilon>0$  such that $\epsilon K\subseteq [-1,1]$. Thus, by \eqref{BasicComparison}, on the event $E_t(\mu)$ \eqref{CouplingEvent}, if $x<y$ and $x,y\in K$, then 
\begin{equation}\label{BlowUpKPZ}
\fh_t^{-\mu}(\epsilon y)-\fh_t^{-\mu}(\epsilon x)\leq \fh_t(\epsilon  y;\fh)-\fh_t(\epsilon  x;\fh)\leq\fh_t^{\mu}(\epsilon  y)-\fh_t^{\mu}(\epsilon x)\,.
\end{equation}
Denote the modulus of continuity of a function $\ff$ by 
$$\omega(\ff,\delta):=\sup_{x, y\in K\,,\,x\neq y\,,\,|x-y|\leq \delta}|\ff(x)-\ff(y)|\,.$$
By \eqref{BlowUpKPZ}, on the event $E_t(\mu)$,
\begin{equation}\label{ModCont}
\omega\left(S_{\sqrt{\epsilon}}\Delta\fh_t,\delta\right)\leq \max\left\{\omega\left(S_{\sqrt{\epsilon}}\Delta\fh^{-\mu}_t,\delta\right)\,,\,\omega\left(S_{\sqrt{\epsilon}}\Delta\fh^{\mu}_t,\delta\right)\right\}\,.
\end{equation}
We note that, for every $\mu\in\R$, 
\begin{equation}\label{BlowUpBM}
S_{\sqrt{\epsilon}}\Delta\fh^{\mu}_t(x)\stackrel{dist.}{=}\mu\epsilon^{1/2}x +\fb(x)\,\mbox{ (as process in $x\in\R$) }\,,
\end{equation}  
and we want to tune $\mu=\mu_\epsilon$ in order to have 
$$\P\left(E_t(\mu)^c\right)\to 0\,\mbox{ and }\,\mu\epsilon^{1/2}\to 0\,,\mbox{ as }\epsilon\to 0\,.$$
By choosing $\mu_\epsilon=\epsilon^{-1/4}$ we have both (using \eqref{BMArgMax} and Lemma \ref{KPZLocalization1} as in \eqref{ControlComparison}), and by \eqref{ModCont} and \eqref{BlowUpBM}, for every $\eta>0$,
$$\P\left(\omega\left(S_{\sqrt{\epsilon}}\Delta\fh_t,\delta\right)>\eta\right)\leq 2\P\left(\omega\left(\fb,\delta\right)>\eta-\delta\epsilon^{1/4}\right)+\P\left(E_t(\mu_\epsilon)^c\right)\,.$$
This shows that for every $\eta_1,\eta_2>0$ there exist $\delta>0$ and $\epsilon_0>0$ such that 
$$\P\left(\omega\left(S_{\sqrt{\epsilon}}\Delta\fh_t,\delta\right)>\eta_1\right)<\eta_2\,,\,\forall\,\epsilon<\epsilon_0\,.$$
Since $S_{\sqrt{\epsilon}}\Delta\fh_t(0)=0$, this implies that the sequence of probability measures in $\cC(K)$ induced by $S_{\sqrt{\epsilon}}\Delta\fh_t$ is tight. On the other hand, by picking $x=0$ in \eqref{BlowUpKPZ}, $\mu_\epsilon=\epsilon^{-1/4}$ and then using \eqref{BlowUpBM}, we see that the finite dimensional distributions of $S_{\sqrt{\epsilon}}\Delta\fh_t$ are converging, as $\epsilon\to 0$, to those of $\fb$, which finishes the proof of \eqref{EqBrown}.

\hfill$\Box$\\

\subsection{Proof of Theorem \ref{AirySheet}} 
Recall that 
$$\cA(x,y)=\cL(x,y)+(x-y)^2\,,$$
where $\cL(x,y):=\cL(x,0;y,1)$, and it is sufficient to prove the analog result for $\cL$. Since $\Delta\cL(0,0)=0$, to prove tightness we only need to control the modulus of continuity of the two-dimensional scalar field $\cL$. Now we can write
\begin{eqnarray*}
\epsilon^{-1/2}\left(\cL(\epsilon  x_2,\epsilon y_2)-\cL(\epsilon  x_1,\epsilon y_1)\right)&=&\epsilon^{-1/2}\left(\cL(\epsilon  x_2,\epsilon y_2)-\cL_1(\epsilon  x_2,\epsilon y_1)\right)\\
&+&\epsilon^{-1/2}\left(\cL(\epsilon  x_2,\epsilon y_1)-\cL(\epsilon  x_1,\epsilon y_1)\right)\,.
\end{eqnarray*}
By the symmetry $\left\{\cL(x,y)\right\}_{(x,y)\in\R^2}\stackrel{dist.}{=}\left\{\cL(y,x)\right\}_{(x,y)\in\R^2}$, it is sufficient to control the supremum of 
$$\epsilon^{-1/2}\left(\cL(\epsilon x,\epsilon  y_2)-\cL(\epsilon x,\epsilon  y_1)\right)\,,$$
over all $(y_1,x),(y_2,x)\in K$ with $|y_1-y_2|\leq \delta$, where $K$ is a fixed compact subset of $\R^2$. Recall that the directed landscape can be expressed as
$$\cL(x,y)=\fh_1(y;\fd_x)\,,\mbox{ where }\,\fd_x(z)=\left\{\begin{array}{ll}
0 & \mbox{ for } z=x\\
-\infty &\mbox{ for } z\neq x\,.\end{array}\right. 
$$ 
Notice also that $Z_1(y;\fd_x)=x$ for all $y\in\R$. Given $K\subseteq \R^2$ compact there exists $\epsilon_0$ such that $\epsilon|x|,\epsilon |y|\leq 1$ for all $(x,y)\in K$ and for all $\epsilon<\epsilon_0$. Hence
$$|Z_1(\epsilon y;\fd_{\epsilon x})|=\epsilon|x|\leq 1\,,\mbox{ for all }(x,y)\in K\,,$$
and, on the event that  
\begin{equation}\label{SheetComp}
Z_t^{-\mu}(1)\leq -1 < 1\leq Z_t^\mu(-1)\,, 
\end{equation}
(as in \eqref{BlowUpKPZ}) we have that for all $(x,y_1)\in K$ and $(x,y_2)\in K$, with $y_1<y_2$,
$$\epsilon^{-1/2}\Big(\fh_1^{-\mu}(\epsilon y_2)-\fh_1^{-\mu}(\epsilon y_1)\Big)\leq \epsilon^{-1/2}\Big(\fh_1(\epsilon  y_2;\fd_{\epsilon x})-\fh_1(\epsilon  y_1;\fd_{\epsilon x})\Big)\leq\epsilon^{-1/2}\Big(\fh_1^{\mu}(\epsilon  y_2)-\fh_1^{\mu}(\epsilon y_1)\Big)\,.$$
For $\mu=\mu_\epsilon=\epsilon^{-1/4}$, \eqref{SheetComp} occurs with high probability as $\epsilon\to 0$, and under \eqref{SheetComp}, for all $x\in \R$ such that $(x,y)\in K$ for some $y\in \R$, we have that  
\begin{equation*}
\omega\left(S_{\sqrt{\epsilon}}\Delta\fh_t(\cdot;\fd_{\epsilon x}),\delta\right)\leq \max\left\{\omega\left(S_{\sqrt{\epsilon}}\Delta\fh^{-\mu}_t,\delta\right)\,,\,\omega\left(S_{\sqrt{\epsilon}}\Delta\fh^{\mu}_t,\delta\right)\right\}\,.
\end{equation*}
From here one can follow the proof of Theorem \ref{Holder_1} to conclude tightness and  marginal local Brownian behaviour. From the same argument, one can get $1/2-$ Holder regularity of the Airy Sheet.  
\newline

To prove independence we have to change the comparison set up, and we do it by splitting the space-time directed landscape at time $s=1/2$. For $x,y\in\R$ consider 
$$Z_{1/2}(x,y)=\cP_{x,0}^{y,1}(1/2)\,,$$ 
i.e. the location at time $s=1/2$ of the rightmost geodesic between $(x,0)$ and $(y,1)$. Thus,  by metric composition \eqref{MetComp}, 
$$\cL(x,y)=\max_{z\in\R}\left\{\cL(x,0;z,1/2)+\cL(z,1/2;y,1)\right\}=\cL(x,0;Z_{1/2},1/2)+\cL(Z_{1/2},1/2;y,1)\,.$$
As in the proof of (iv)-Proposition \ref{PropBasic}, we have monotonicity of geodesics as follows: for all $x_1\leq x_2$ and $y_1\leq y_2$ then
$$\cP_{x_1,0}^{y_1,1}(s)\leq\cP_{x_2,0}^{y_2,1}(s)\,,\,\forall\,s\in[0,1]\,,$$
and, in particular, 
\begin{equation}\label{GeoMon}
Z_{1/2}(x_1,y_1)\leq Z_{1/2}(x_2,y_2)\,.
\end{equation}
Let
$$\bar\fh_{1/2+}(y;\fh):=\max_{z\in\R}\left\{\fh(z)+\cL(z,1/2;y,1)\right\}\,\mbox{ and }\,\bar\fh_{1/2-}(x;\fh):=\max_{z\in\R}\left\{\fh(z)+\cL(x,0;z,1/2)\right\}\,.$$
Then, by metric composition \eqref{MetComp}, 
$$\cL(x,y)= \bar\fh_{1/2+}(y;\bar\fh_{+,x})\,\mbox{ and }\,\cL(x,0)= \bar\fh_{1/2-}(x;\bar\fh_-)\,,$$
where 
$$\bar\fh_{+,x}(z)=\cL(x,0;z,1/2)=\fh_{1/2}(z;\fd_x)\,\mbox{ and }\,\bar\fh_{-}(z)=\cL(z,1/2;0,1)\,.$$ 
Therefore,
$$\cL(x,y)-\cL(0,0)=\cL(x,y)-\cL(x,0)+\cL(x,0)-\cL(0,0)=\Delta\bar\fh_{1/2+}(y;\bar\fh_{+,x})+\Delta\bar\fh_{1/2-}(x;\bar\fh_{-})\,.$$
The trick now is to pick $\fb_1$ and $\fb_2$, two independent copies of $\fb$, and then apply the coupling method to compare simultaneously $\Delta\bar\fh_{1/2+}(y;\bar\fh_{+,x})$ with $\Delta\bar\fh_{1/2+}(y;\fb^\mu_1)$, and $\Delta\bar\fh_{1/2-}(y;\bar\fh_{-})$ with $\Delta\bar\fh_{1/2-}(y;\fb^\mu_2)$. By time independence and stationarity \eqref{DirSheet} of the directed landscape, we clearly have that $\bar\fh_{1/2+}(\cdot;\fb^\mu_1)$ and $\bar\fh_{1/2-}(\cdot;\fb^\mu_2)$ are independent processes, and 
$$\Delta\bar\fh_{1/2+}(\cdot;\fb^\mu_1)\stackrel{dist.}{=}\fb^\mu\stackrel{dist.}{=}\Delta\bar\fh_{1/2-}(\cdot;\fb^\mu_2)\,.$$    
Let 
$$\bar Z_{1/2+}(y,\fh):=\max\argmax_{z\in\R}\left\{\fh(z)+\cL(z,1/2;y,1)\right\}\,,$$
and 
$$\bar Z_{1/2-}(x,\fh):=\max\argmax_{z\in\R}\left\{\fh(z)+\cL(x,0;z,1/2)\right\}\,.$$
Hence, 
\begin{equation*}
\bar Z_{1/2+}(y,\bar\fh_{+,x})=Z_{1/2}(x,y)\,\mbox{ and }\,\bar Z_{1/2-}(x,\bar\fh_{-})=Z_{1/2}(x,0)\,.
\end{equation*}
Let
$$\bar E_{1/2+}(\mu):=\left\{\bar Z_{1/2+}(-1,\fb^{\mu}_1)\geq Z_{1/2}(1,1)\,\mbox{ and }\, \bar Z_{1/2+}(1,\fb^{-\mu}_1)\leq Z_{1/2}(-1,-1) \right\}\,,$$ 
and 
$$\bar E_{1/2-}(\mu):=\left\{\bar Z_{1/2-}(-1,\fb^{\mu}_2)\geq Z_{1/2}(1,0)\,\mbox{ and }\, \bar Z_{1/2+}(1,\fb^{-\mu}_2)\leq Z_{1/2}(-1,0) \right\}\,.$$ 
For a compact set $K\subseteq \R^2$ we can chose again $\epsilon_0$ so that $\epsilon|x|,\epsilon |y|\leq 1$ for all $(x,y)\in K$ and for all $\epsilon<\epsilon_0$. Thus, by \eqref{GeoMon},
$$Z_{1/2}(-1,-1)\leq Z_{1/2}(\epsilon x,\epsilon y)\leq Z_{1/2}(1,1)\,\,\mbox{ and }\,\,Z_{1/2}(-1,0)\leq Z_{1/2}(\epsilon x,0)\leq Z_{1/2}(1,0)\,\,.$$
Denote $\bar\fh^{\pm\mu}_{+1/2}(\cdot)\equiv\bar\fh_{+1/2}(\cdot;\fb^{\pm\mu}_1)$ and  $\bar\fh^{\pm\mu}_{-1/2}(\cdot)\equiv\bar\fh_{-1/2}(\cdot;\fb^{\pm\mu}_2)$. On the event $\bar E_{1/2+}(\mu)$, for all $(x,y)\in K$, if $0<y$ then 
$$\bar\fh_{1/2+}^{-\mu}(\epsilon y)-\bar\fh_{1/2+}^{-\mu}(0)\leq\bar\fh_{1/2+}(\epsilon  y;\fd_{\epsilon x})-\bar\fh_{1/2+}(0;\fd_{\epsilon x})\leq\bar\fh_{1/2+}^{\mu}(\epsilon  y)-\bar\fh_{1/2+}^{\mu}(0)\,,$$
while if $y<0$ then
$$\bar\fh_{1/2+}^{\mu}(\epsilon y)-\bar\fh_{1/2+}^{\mu}(0)\leq \bar\fh_{1/2+}(\epsilon  y;\fd_{\epsilon x})-\bar\fh_{1/2+}(0;\fd_{\epsilon x})\leq\bar\fh_{1/2+}^{-\mu}(\epsilon  y)-\bar\fh_{1/2+}^{-\mu}(0)\,.$$
On the event $\bar E_{1/2-}(\mu)$ for all $(x,y)\in K$, if $0<y$ then 
$$\bar\fh_{1/2-}^{-\mu}(\epsilon y)-\bar\fh_{1/2-}^{-\mu}(0)\leq\bar\fh_{1/2-}(\epsilon  y;\bar\fh_-)-\bar\fh_{1/2-}(0;\bar\fh_-)\leq\bar\fh_{1/2-}^{\mu}(\epsilon  y)-\bar\fh_{1/2-}^{\mu}(0)\,,$$
while if $y<0$ then
$$\bar\fh_{1/2-}^{\mu}(\epsilon y)-\bar\fh_{1/2-}^{\mu}(0)\leq \bar\fh_{1/2-}(\epsilon  y;\bar\fh_-)-\bar\fh_{1/2-}(0;\bar\fh_-)\leq\bar\fh_{1/2-}^{-\mu}(\epsilon  y)-\bar\fh_{1/2-}^{-\mu}(0)\,.$$
Thus, for $\mu=\mu_\epsilon=\epsilon^{-1/4}$, on the event $\bar E_{1/2+}(\mu)\cap\bar E_{1/2-}(\mu)$, one can approximate the finite dimensional distributions of $(S_{\sqrt{\epsilon}}\Delta\bar\fh_{1/2+},S_{\sqrt{\epsilon}}\Delta\bar\fh_{1/2-})$ using the finite dimensional distributions of $(\fb_1,\fb_2)$ (as in the proof of Theorem \ref{Holder_1}). Since
$$\P\left(\bar E_{1/2+}(\mu)\cap\bar E_{1/2-}(\mu)\right)\to 1\mbox{ as }\epsilon\to 0\,,$$
this finishes the proof Theorem \ref{AirySheet}.

\subsection{Proof of Theorem \ref{Coupling}} 
Recall  \eqref{DriftBM} and \eqref{DriftBM_1}. By Proposition \ref{Attractiveness} (attractiveness), 
$$\Delta\fh_t^{-\mu}(x)\leq \Delta \fh_t^{0}(x)\leq \Delta\fh_t^{+\mu}(x)\,,\mbox{ for $x\geq 0$}\,,$$
and
$$\Delta\fh_t^{+\mu}(x)\leq \Delta\fh_t^{0}(x)\leq \Delta\fh_t^{-\mu}(x)\,,\mbox{ for $x\leq 0$}\,.$$
Furthermore \footnote{It also follows from Proposition \ref{Attractiveness} that $\Delta\fh_t^{+\mu}(x)-\Delta\fh_t^{-\mu}(x)$ is a nondecreasing function of $x$.},
$$0\leq\Delta\fh_t^{+\mu}(x)-\Delta\fh_t^{-\mu}(x)\leq \Delta\fh_t^{+\mu}(a)-\Delta\fh_t^{-\mu}(a)\,,\,\forall x\in[0,a]\,,$$
and 
$$0\leq\Delta\fh_t^{-\mu}(x)-\Delta\fh_t^{+\mu}(x)\leq \Delta\fh_t^{-\mu}(-a)-\Delta\fh_t^{+\mu}(-a)\,,\,\forall x\in[-a,0]\,.$$
By time invariance \eqref{stat_1}, $\Delta\fh_t^{\pm\mu}(x)$ is a two-sided Brownian motion with drift $\pm\mu$. Hence
$$\E\left(\Delta\fh_t^{+\mu}(a)-\Delta\fh_t^{-\mu}(a)\right)=\E\left(\Delta\fh_t^{-\mu}(-a)-\Delta\fh_t^{+\mu}(-a)\right)= 2\mu a\,.$$

Consider the event $E_t(\mu)$ \eqref{CouplingEvent}. By \eqref{BasicComparison}, 
$$\Delta\fh_t^{-\mu}(x)\leq \Delta \fh_t(x;\fh)\leq \Delta\fh_t^{+\mu}(x)\,,\mbox{ for $x\in [0,a]$}\,,$$
and
$$\Delta\fh_t^{+\mu}(x)\leq \Delta\fh_t(x;\fh)\leq \Delta\fh_t^{-\mu}(x)\,,\mbox{ for $x\in [-a, 0]$}\,.$$
\newline
Thus, if $E_t(\mu)$ occurs then both $\Delta\fh^{0}_t(\cdot)=\Delta\fh_t(\cdot;\fb)$ and $\Delta\fh_t(\cdot;\fh)$ are sandwiched by $\Delta\fh^{\pm\mu}_t(\cdot)$, which implies the following uniform control on the distance between $\Delta\fh_t(\cdot;\fh)$ and $\Delta\fh_t(\cdot;\fb)$:
$$\sup_{x\in[-a,a]}|\Delta\fh_{t}(x;\fh)-\Delta \fh_{t}(x;\fb)|\leq I_t(a)\,,$$  
where 
$$0\leq I_t(a)=\Delta\fh_t^{+\mu}(a)-\Delta\fh_t^{-\mu}(a)+\Delta\fh_t^{-\mu}(-a)-\Delta\fh_t^{+\mu}(-a)\,.$$
Therefore, using Markov inequality and that $\E\left(I_t(a)\right)=4\mu a$, 
$$\P\left( \sup_{x\in[-a,a]}|\Delta\fh_{t}(x;\fh)-\Delta \fh_{t}(x;\fb)|>\eta\sqrt{a}\right)\leq \P\left(E_t(\mu)^c\right)+\frac{\E\left(I_t(a)\right)}{\eta\sqrt{a}}=\P\left(E_t(\mu)^c\right)+\frac{4\mu\sqrt{a}}{\eta}\,.$$
\newline

In order to make this inequality useful, we have to chose $\mu=\mu_t$ in such way that 
$$\P\left(E_t(\mu)^c\right)\to 0\,\,\mbox{ and }\,\,\mu\sqrt{a}\to 0\,,\,\mbox{ as }t\to\infty\,$$ 
(we allow $a=a_t$ as well). For $t\geq a^{3/2}$ we have that $\pm at^{-2/3}$ does not play any rule in the asymptotic analysis of $E_t(\mu)$ (recall \eqref{BMArgMax}). By Lemma \ref{KPZLocalization2}, we know that 
$$\P\left(|Z_t(\pm a;\fh)|>rt^{2/3}\right)\to 0\,,\mbox{ as }r\to\infty\,,$$ 
(uniformly in $t$). Thus, by \eqref{BMArgMax}, $E_t(\mu)$ should occur with high probability, as soon as  $\pm \mu t^{1/3}\to\pm\infty$. By setting $\mu= r (4t^{1/3})^{-1}$, for some $r=r_t\to\infty$, then   
$$4\mu\sqrt{a}=r(at^{-2/3})^{1/2}\,.$$
A natural choice is $r_t=(at^{-2/3})^{-\delta}$ with $\delta\in(0,1/2)$, and for the sake of simplicity we take $\delta=1/4$, which yields to
\begin{equation}\label{IneMain}
\P\Big( \sup_{x\in[-a,a]}|\Delta\fh_{t}(x;\fh)-\Delta \fh_{t}(x;\fb)|>\eta\sqrt{a}\Big)\leq \P\left(E_t(\mu)^c\right)+\frac{1}{\eta r_t}\,.
\end{equation}
Therefore, Theorem \ref{Coupling} is a consequence of \eqref{IneMain} and Lemma \ref{CompEvent} below. 

\begin{lem}\label{CompEvent}
Let $\mu :=r(4t^{1/3})^{-1}$. Then, under assumption \eqref{Assump1}, there exists a function $\phi$, that does not depend on $a,t>0$, such that for all $t\geq \max\{c^3,a^{3/2}\}$ 
$$\P\left(E_t(\mu)^c\right)\leq \phi(r)\,\,\mbox{ and }\,\,\lim_{r\to\infty}\phi(r)=0\,.$$
\end{lem}

\noindent\paragraph{\bf Proof}
By the definition of $E_t(\mu)$, 
$$E_t(\mu)^c\cap\left\{|Z_t(a;\fh)|\leq \frac{r}{16}t^{2/3}\right\}\subseteq\left\{Z^{\mu}_t(-a)\leq \frac{r}{16}t^{2/3}\right\}\cup\left\{Z^{-\mu}_t(a)\geq -\frac{r}{16}t^{2/3}\right\}\,,$$ 
and hence,
\begin{equation}\label{CompEvent1}
\P\left(E_t(\mu)^c\right)\leq \P\left(|Z_t(a;\fh)|>\frac{r}{16}t^{2/3}\right)+\P\left(Z^{\mu}_t(-a)\leq \frac{r}{16}t^{2/3}\right)+\P\left(Z^{-\mu}_t(a)\geq -\frac{r}{16}t^{2/3}\right)\,.
\end{equation}
By Lemma \ref{KPZLocalization2}, we only need to show that there exists a function $\phi_2$, that does not depend on $a,t>0$, such that for all $t\geq \max\{c^3,a^{3/2}\}$ 
$$\max\left\{\P\left(Z^{\mu}_t(-a)\leq \frac{r}{16}t^{2/3}\right)\,,\,\P\left(Z^{-\mu}_t(a)\geq -\frac{r}{16}t^{2/3}\right)\right\}\leq \phi_2(r)\,,$$
and $\lim_{r\to\infty}\phi_2(r)=0$. Since $\mu :=r(4t^{1/3})^{-1}$ and $t\geq a^{3/2}$, by \eqref{BMArgMax}, 
$$\P\left(Z^{\mu}_t(-a)\leq \frac{r}{16}t^{2/3}\right)=\P\left(Z^{0}_1(0)\leq -\frac{r}{16}+at^{-2/3}\right)\leq\P\left(Z^{0}_1(0)\leq -\left(\frac{r}{16}-1\right)\right) \,,$$
and 
$$\P\left(Z^{-\mu}_t(a)\geq -\frac{r}{16}t^{2/3}\right)=\P\left(Z^{0}_1(0)\geq \frac{r}{16}-at^{-2/3}\right)\leq \P\left(Z^{0}_1(0)\geq \frac{r}{16}-1\right)\,,$$
which allows us to take $\phi_2(r):=\P\left(|Z^{0}_1(0)|>\frac{r}{16}-1\right)$. Therefore, together with \eqref{CompEvent1}, this shows that 
$$\nonumber\P\left(E_t(\mu)^c\right)\leq \phi_1(r)+2\phi_2(r)\,.$$

\hfill$\Box$\\


\begin{thebibliography}{10}

\bibitem{AmCoQu}
\textsc{G. Amir, I. Corwin and J. Quastel.} Probability distribution of the free energy of the continuum directed random polymer in 1 + 1 dimensions. \newblock\emph{Commun. Pure Appl. Math.} {\bf 64}:466--537 (2011).

\bibitem{BaDeJo}
\textsc{J. Baik, P.~A Deift and K. Johansson}. On the distribution of the length of the longest increasing subsequence of random permutations.
\newblock \emph{J. Amer. Math. Soc.} {\bf 12}:1119--1178 (1999).

\bibitem{BaCaSe}
\textsc{M. Bal\'azs, E.~A. Cator and T. Sepp\"al\"ainen}. Cube root fluctuations for the corner growth model associated to the exclusion process.
\newblock \emph{Elect. J. Probab.} {\bf 11}:1094--1132 (2006).

\bibitem{BaQuSe} 
\textsc{M. Bal\'azs, J. Quastel, T. Sepp\"al\"ainen}. Scaling exponent for the Hopf-Cole solution of KPZ/Stochastic Burgers. 
\newblock \emph{J. Amer. Math. Soc.} {\bf 24}:683--708 (2011).

\bibitem{BaSe}
\textsc{M. Bal\'azs and T. Sepp\"al\"ainen.} Fluctuation bounds for the asymmetric simple exclusion process.
\newblock \emph{ALEA Lat. Am. J. Probab. Math. Stat.} {\bf 6}:1--24, (2009).

\bibitem{BoFePrSa}
\textsc{A. Borodin, P.~L. Ferrari, M. Pr\"ahofer and T. Sasamoto}. Fluctuation properties of the TASEP with periodic initial configuration. 
\newblock \emph{J. Statist. Phys.} {\bf 129}:1055--1080 (2007).

\bibitem{CaGr}
\textsc{E.~A Cator and P. Groeneboom}. Second class particles the cube root asymptotics  for Hammersley's process. \newblock\emph{Ann. Probab.} {\bf 34}:1273--1295 (2006).

\bibitem{CaPi}
\textsc{E.~A. Cator and L.~P.~R. Pimentel}. On the local fluctuations of last-passage percolation models. \newblock\emph{Stoch. Proc. Appl.} {\bf 125}:879--903 (2012).


\bibitem{Co}
\textsc{I. Corwin}. The Kardar-Parisi-Zhang equation and universality class.  \newblock\emph{Random Matrices Theory
Appl.} {\bf 1}(1):1130001, 76, 2012.

\bibitem{CoHa}
\textsc{I. Corwin and A. Hammond}. Brownian Gibbs property for Airy line ensembles. \newblock\emph{Invent. Math.} {\bf 195}:441--508 (2014).


\bibitem{CoLiWa}
\textsc{I. Corwin, Z. Liu and D. Wang}. Fluctuations of TASEP and LPP with general initial data. \newblock\emph{Ann. Appl. Probab.} {\bf 26}:2030--2082 (2016).


\bibitem{CoQuRe}
\textsc{I. Corwin, J. Quastel and D. Remenik}. Renormalization fixed point of the KPZ universality class. \newblock\emph{J. Stat. Phys.} {\bf 160}:815--834  (2015).

\bibitem{DaOrVi}
\textsc{D. Dauvergne, J. Ortmann and B\'alint Vir\'ag}. The directed landscape. \newblock Available from arXiv:1812.00309.

\bibitem{FeOc}
\textsc{P.~L. Ferrari and A. Occelli}. Universality of the GOE Tracy-Widom distribution for TASEP with arbitrary particle density. \newblock \emph{Elect. J. Probab.} {\bf 23}, no 51:1--24 (2018).

\bibitem{Ha}
\textsc{T.~E. Harris}. Additive set-valued Markov processes and graphical methods.\newblock \emph{Ann. Probab.} {\bf 6}:355--378 (1978).

\bibitem{Jo1}
\textsc{K. Johansson}. Shape fluctuations and random matrices. \newblock \emph{Comm. Math. Phys.} {\bf 209}:437--476 (2000).

\bibitem{Jo2}
\textsc{K. Johansson}. Discrete Polynuclear Growth and Determinantal processes. \newblock \emph{Comm. Math. Phys.} {\bf 242}:277--239 (2003).

\bibitem{KPZ}
\textsc{M. Kardar, G. Parisi, Y.~-C. Zhang}. Dynamic scaling of growing interfaces. \newblock\emph{Phys. Rev. Lett.} {\bf 56}:889--892 (1986).

\bibitem{MaQuRe}
\textsc{K. Matetski, J. Quastel, and D. Remenik}. The KPZ fixed point.
\newblock Available from arXiv:1701.00018.

\bibitem{Pi1}
\textsc{ L.~P.~R. Pimentel}. On the location of the maximum of a continuous stochastic process.
\newblock \newblock\emph{J. Appl. Probab.} {\bf 173}: 152--161 (2014).

\bibitem{Pi}
\textsc{ L.~P.~R. Pimentel}. Local behavior of Airy processes.
\newblock \newblock\emph{J. Stat. Phys.} {\bf 173}: 1614--1638 (2018).

\bibitem{Pi2}
\textsc{ L.~P.~R. Pimentel}. Ergodicity of  the KPZ fixed point.
\newblock \newblock Available from arXiv:1708.06006.

\bibitem{PrSp}
\textsc{M. Pr\"ahofer and H. Spohn}. Scale invariance of the PNG droplet and the Airy process. \newblock \emph{J. Stat. Phys.} {\bf 108}:1071--1106 (2002).

\bibitem{Se}
\textsc{T. Sepp\"al\"ainen}. Scaling for a one-dimensional directed polymer with boundary conditions. \newblock\emph{Ann. Probab.} {\bf 40} 19--73 (2012).

\bibitem{Ta}
\textsc{K.~A. Takeuchi, M. Sano, T. Sasamoto and  H. Spohn}. Growing interfaces uncover universal fluctuations behind scale invariance. \newblock \emph{Scientific Reports (Nature)}. {\bf 1}: 34 (2011).


\end{thebibliography}
\end{document}